\documentclass[12pt]{article}
\usepackage{amssymb}
\usepackage{myart}

\normalbaselines
\input tcilatex

\begin{document}

\title{New crisis in geometry? }
\author{Yuri A.Rylov}
\date{Institute for Problems in Mechanics, Russian Academy of Sciences,\\
101-1, Vernadskii Ave., Moscow, 119526, Russia.\\
e-mail: rylov@ipmnet.ru\\
Web site: {$http://rsfq1.physics.sunysb.edu/\symbol{126}rylov/yrylov.htm$}\\
or mirror Web site: {$http://gasdyn-ipm.ipmnet.ru/\symbol{126}%
rylov/yrylov.htm$}}
\maketitle

\begin{abstract}
The first crisis in the geometry arose in the beginning of XIXth century,
when the mathematicians rejected the non-Euclidean geometry as a possible
geometry of the real world. Now we observe unreasonable rejection of the
non-Riemannian geometry by the official representatives of the contemporary
geometry. Class of the Riemannian geometries appears to be too narrow for
physical applications. The microcosm physics needs expansion of the class of
possible geometries appropriate for the role of space-time geometry. In the
framework of the non-Riemannian geometry one can construct the space-time
geometry, where the motion of free particles is primordially stochastic, and
this stochasticity depends on the particle mass. At the same time the
geometry in itself is not stochastic in the sense that the space-time
intervals are deterministic. Principles of quantum mechanics can be deduced
from such a space-time geometry. The crisis situation in geometry appears to
be connected with some preconceptions concerning the foundation of the
geometry. The preconceptions as well as the crisis generated by them are not
purely scientific phenomena. The human factor (social aspect) is rather
strong in the crisis phenomena. The preconceptions and the human factor
appear to be so strong, that usual logical arguments are not perceived, and
the usual formal mathematical language appears to be inappropriate for
perception of an analysis of the crisis origin and of a possibility of its
overcoming. In the paper the history and motives of the non-Riemannian
geometry construction are presented. There is a hope that such a less formal
way of presentation helps to understand and to overcome the existing
preconceptions.
\end{abstract}

\section{Introduction}

The geometry as well as the mathematics as a whole is an exact logical
science. It seems that all problems in geometry can be resolved by means of
logical reasonings, and crisis phenomena seem to be impossible in geometry.
Speaking on crisis phenomena, I mean such a situation, when the most of
geometricians do not want to accept the well founded statements of geometry,
having far reached consequences for the microcosm physics. Some researchers
wrote on crisis phenomena in physics and mathematics in the end of XXth
century, keeping in mind such a situation, when the science development
becomes slow. S.P. Novikov \cite{N2000} considers social origin of the
crisis ( insufficient education, unsatisfactory organization of scientific
investigations, etc.). I do consider mainly the scientific origin of the
crisis. I keep in mind such a situation, which took place in the XIXth
century, when the mathematician rejected to accept the non-Euclidean
geometry as a possible geometry of the real world without any reasonable
arguments. The subsequent development of the non-Euclidean geometry and its
application in the theory of gravitation showed that the rejection of the
non-Euclidean geometry was a mistake. But it was not simply a mistake, it
was a crisis in geometry, because this mistake was overall. The viewpoint of
advocates of the non-Euclidean geometry was blamed by all (or almost all)
mathematicians.

Now we have the same situation in the non-Riemannian geometry. Why do I
think, that the mathematicians do not accept the conception of
non-Riemannian geometry? The conceptual paper on foundation of the
non-Riemannian geometry was sent in turn to several leading mathematical
journals. It was rejected without any motivation. To be not naked I present
correspondence concerning submission of this paper to one of such journals.
The name of the mathematical journal and those of the editor are coded.
Personality of editors is of no importance, as far as we discuss the
phenomena, characteristic for the whole geometrical community.

\qquad \qquad \qquad \qquad \qquad \qquad \qquad \qquad \qquad \qquad \qquad
\qquad \qquad \textit{June 21, 2004}

\textit{Dear Prof. S.,}

\textit{I should like to submit my paper "Coordinateless description and
deformation principle as foundations of physical geometry" for publication
in JABCD.}

\textit{I am a physicist (not mathematician) dealing with the space-time and
its geometry, and geometry is interesting for me as a tool for my
investigations in physics. I was told that such a problem as foundation of
geometry is considered to be solved, and none of mathematicians deals with
it. When I have submit my paper to the Editor-in-Chief, he recommended me to
submit the paper to you.}

\textit{From my (physical) point of view the simpler is geometry, the better
is its application to physical problems. The geometry foundation appears to
be very simple (may be, indecently simple), because I do not try to repeat
Euclidean constructions at other axiomatics. Instead, I use ready Euclidean
geometrical construction, deforming them in proper way. Any new way of
deformation generates a new physical geometry. As a result one succeeds to
construct such geometries, which cannot be constructed by conventional
methods. It is such geometries that appear to be interesting as space-time
geometries. For instance, one can construct such a deterministic space-time
geometry, where motion of free particles is stochastic and the particle mass
is geometrized.}

\textit{I am understanding that the too simple method of the geometry
construction may be not interesting for mathematicians, but it is very
interesting from the viewpoint of applications to physics.}

\textit{Sincerely yours,}

\textit{Yuri Rylov}

\textit{P.S. Confirm reception, please.}

\bigskip

\textit{\qquad \qquad \qquad \qquad \qquad \qquad \qquad \qquad \qquad
\qquad \qquad \qquad \qquad September 30, 2004}

\textit{Dear Professor Rylov:}

\textit{This is to acknowledge that we have received paper submitted to
JABCD entitled "Coordinateless description and deformation principle as
foundations of physical geometry." Please keep me informed of any change in
your address (permanent or temporary).}

\textit{Sincerely,}

\textit{S., Editor, JABCD}

\bigskip

\textit{Professor Yuri A. Rylov}

\textit{Institute for Problems in Mechanics,}

\textit{\ Russian Academy of Sciences}

\textit{101--1, Vernadskii Avenue Moscow, 119526, Russia}

\textit{email:rylov@ipmnet.ru}

\textit{\ \qquad \qquad \qquad \qquad \qquad \qquad \qquad \qquad \qquad
\qquad \qquad \qquad \qquad October 1, 2004}

\textit{Dear Professor Rylov:}

\textit{There was a misunderstanding your paper has been looked at, and I am
sorry to inform you that I have received a quick evaluation of your paper
\textquotedblleft Coordinateless description and deformation principle as
foundations of physical geometry,\textquotedblright\ indicating that your
work is possibly interesting and worthy of publication in a specialist
journal. However, it does not quite reach the high standards of JABCD for
which we expect results with a more profound impact and with rather wider
mathematical implications.}

\textit{Under the circumstances, I will not initiate a formal refereeing
process because the chances that your paper will be accepted in JABCD are
very slim.}

\textit{I hope that you will find an alterative way for its diffusion.}

\textit{Sincerely,}

\textit{S.}

\textit{\ \qquad \qquad \qquad \qquad \qquad \qquad \qquad \qquad \qquad
\qquad \qquad \qquad \qquad October 2, 2004}

\textit{Dear Professor S.,}

\textit{Thank you for your kind letter, concerning my submitted paper
"Coordinateless description and deformation principle as foundations of
physical geometry". Unfortunately, your letter is too polite, and I have not
understood, whether my paper has been rejected, or not yet.}

\textit{If my paper has been rejected, I should like to know motives of the
rejection.}

\textit{You wrote that my paper does not reach the high standards of JABCD.
But it is a naked proofless statement.}

\textit{I am understanding, that it is very difficult to find a qualified
referee for a review of my paper, because my paper investigate geometry in
itself, whereas using coordinate methods, most geometricians investigate
only the methods of the geometry description, but not the geometry in
itself. Of course, it is difficult for them to review the paper on geometry.}

\textit{Coordinateless method of the geometry description is more effective,
than the conventional coordinate method of description. Euclid created his
Euclidean geometry, using the coordinateless method. Coordinateless method
allows one to create geometry with intransitive property of parallelism. The
real space-time geometry is a geometry with the intransitive parallelism
property, and such geometry cannot be constructed in the framework of the
conventional conception of Riemannian geometry.}

\textit{Physics needs geometries, which could be used as space-time
geometries. Geometricians must create such geometries. Unfortunately, they
could not create geometries appropriate for the correct description of the
space-time. When physicists create the class of geometries appropriate for
the space-time description, it appears that such works "do not reach the
high standards of JABCD" and "rather wider mathematical implications"}

\textit{Once 150 years ago mathematicians earned fame, when they rejected
non-Euclidean geometry. I believe they should not earn such a fame again,
rejecting non-Riemannian geometry.}

\textit{Please, if my paper has been rejected, confirm this fact, and if
possible motivate it.}

\textit{Sincerely yours,}

\textit{Yuri Rylov}

\textit{mailto: rylov@ipmnet.ru}

\textit{\ \qquad \qquad \qquad \qquad \qquad \qquad \qquad \qquad \qquad
\qquad \qquad \qquad \qquad October 5, 2004}

\textit{Dear Dr. Rylov,}

\textit{the rejection does not mean that the paper is without value. It
simply means that we believe it is more appropriate to appear in a journal
specialized on the foundations of geometry.}

\textit{Sincerely,}

\textit{S.}

\textit{\ \qquad \qquad \qquad \qquad \qquad \qquad \qquad \qquad \qquad
\qquad \qquad \qquad \qquad October 5, 2004}

\textit{Dear Professor S.,}

\textit{Thank you for your definite answer concerning my paper
"Coordinateless description and deformation principle as foundations of
physical geometry".}

\textit{I would be very obliged to you, if you could point out a
mathematical journal devoted to foundation of geometry.}

\textit{You see, I am a physicist, and situation in mathematics is unknown
for me. In 2002 I participated in the work of the International meeting on
geometry in Petersburg. I tried to discuss the problem of the geometry
foundation with participants of the conference. Unfortunately, nobody of
mathematicians did not like to discuss such problems. They said that the
problem of the geometry foundation had been solved many years ago, and now
it is interesting for nobody.}

\textit{In such a situation I could not image, that there are mathematical
journals devoted to the geometry foundation. You are quite right, that my
paper should be published in a special journal, devoted to the problem of
the geometry foundation. Unfortunately, I have not known about existence of
such journals, and submit my paper in general mathematical journals.}

\textit{You are the first person, who said me about such journals.}

\textit{Please, point out me at least one of such journals.}

\textit{Thank you in advance.}

\textit{Sincerely yours,}

\textit{Yuri Rylov }

Unfortunately, there was no answer, and I am forced to think that Professor
S. recommended me to submit my paper in defunct journals.

The paper, submitted to JABCD can be found in Archives of LANL \cite{R2004}.
It contains two essential points: (1) it is shown that the Euclidean
geometry can be presented in terms and only in terms of the world function $%
\sigma _{\mathrm{E}}$, (2) any other (in general, non-Riemannian) geometry
is obtained simply by means of a change $\sigma _{\mathrm{E}}\rightarrow
\sigma $, in the presentation of the Euclidean geometry. Here $\sigma $ is
the world function of the geometry in question. Such a construction of the
non-Riemannian geometries is very simple, because it does not need any
logical reasonings. There are no points, against which one could object.
Besides, some difficulties of the conventional Riemannian geometry (the
problem of convexity and lack of absolute parallelism) are absent in the
suggested form of the non-Riemannian geometry. These difficulties speak for
an imperfection of the Riemannian geometry.

Rejection of the paper cannot be explained by a specific editorial policy of
JABCD or by a lack of the scientific scrupulosity of professor S., because
several leading mathematical journals rejected different versions of this
paper without a serious motivation or without any motivation at all. In such
a situation one is forced to speak about the position of the geometrical
community as a whole. This position hinders from a development of the
non-Riemannian geometry in favor of the conventional Riemannian geometry.
Thus, the history of the non-Euclidean geometry development is repeating. It
means a new crisis in the geometry development.

From physical viewpoint a development of the non-Riemannian geometry is
quite necessary. Experiments show that the free microparticle (electrons,
protons, etc.) move stochastically. The motion of free particles is
determined \textit{only} \textit{by the space-time geometry}. In the
beginning of the XXth century, when this problem appeared, there were no
such geometry, where the particle motion be primordially stochastic, and
this stochasticity depend on the particle mass. Then the problem of the
stochastic particle motion was solved by means by introduction of the
additional supposition, known as the quantum mechanics principles.
Investigations show \cite{R03} that the introduction of QM principles is not
necessary. The problem of the primordial stochastic motion of microparticles
can be solved by means of the correct choice of the non-Riemannian
space-time geometry, where the quantum constant is a parameter of the
geometry. Such a conception is more fundamental, than the conventional
quantum theory, because it does not contain QM principles, which are valid
only at the description of the nonrelativistic phenomena.

Mistakes may appear in any scientific conception, and relation to possible
mistakes determines the investigation strategy. There are two different
kinds of the investigation strategies: (1) the Newtonian strategy with its
slogan "Hypotheses non fingo" and (2) the Ptolemaic strategy which does not
burden itself by a search of possible mistakes, compensating them by means
of additional hypotheses. The Ptolemaic strategy is effective at the
investigation of single new phenomena, when one needs to invent new
conceptions for their description. The Ptolemaic strategy is not effective
in the construction of fundamental theory, when a serious systematization of
many phenomena is needed.

Fundamental scientific theories, i.e. theories containing a systematization
of a broad class of phenomena, are created by means of the strategy of Sir
Isaac Newton. The fundamental theory is to contain the minimal number of
fundamental statements. All statements of such a theory are deduced from the
fundamental statements by means of the logic. The logical reasonings may be
long, but they are to be valid. The fundamental theory possesses a great
predictability. If the fundamental theory has serious problems with
prediction and explanation of a series of new phenomena, the Newtonian
strategy recommends one to search for mistake in the fundamental statements
of the theory. After the finding of the mistake it must be corrected. If the
mistake is serious, its discovery may lead to a change of the fundamental
statements of the theory. The search of a possible mistake is very
difficult. If we cannot find the mistake, we are forced to invent a new
hypothesis, which compensates the possible mistake.

When the invention of new hypotheses becomes a system, one should speak on
the Ptolemaic investigation strategy. On the one hand, at the Ptolemaic
strategy there is no necessity to search for possible mistakes in the
foundation of the theory (maybe, such mistakes are absent at all!), and it
is an advantage of this strategy. On the other hand, the predictability
decreases together with the number of additional hypotheses. It is a defect
of the Ptolemaic strategy. Creation of the fundamental theory having the
great predictability is practically impossible by means of the Ptolemaic
strategy. The development of the Ptolemaic conception, containing hidden
mistakes in the fundamental statements, appears to be difficult. It means
that the mistakes in the fundamental statements should be yet discovered and
corrected.

Evaluating a scientific theory, one uses different criteria. Comparing two
Ptolemaic conceptions, based on different systems of fundamental statements,
one takes into account mainly the number of explained phenomena. The number
and compatibility of the fundamental statements is less important. For
instance, if for explanation of a new phenomenon on needs to introduce a new
additional statement, and one succeeds to find this statement, this
explanation is considered to be a serious achievement of the existing
theory. But this viewpoint is only the opinion of the advocates of the
existing theory. From viewpoint of advocates of the strict Newtonian
approach such an explanation is simply a fitting, but not a progress of the
theory.

On the contrary, if one succeeds to find a mistake in the fundamental
statements and reduce the number of them, representing some of previous
fundamental statements as a corollaries of the new system of fundamental
statements, the advocates of the Newtonian approach consider this fact as a
serious achievement of the theory. But the advocates of the Ptolemaic
approach put the question in other form: "What new phenomena are explained
by the new system of fundamental statements?" Advocates of the Ptolemaic
approach can easily answer such questions, because they invent their
additional hypotheses to answer such questions. On the contrary, such
questions present difficulties for advocates of the Newtonian approach,
because starting from the new system of fundamental statements, one needs to
produce long logical reasonings and calculations. For instance, in the
celestial mechanics the change of the Ptolemaic doctrine by the Copernicus
one generates the discovery of the Newtonian gravitation law. It was the
first corollary of the Copernicus doctrine, which cannot be obtained in the
framework of the Ptolemaic doctrine. This discovery meant the triumph of the
Copernicus doctrine, but this triumph took place more than a century later
after creation of the Copernicus doctrine. This example shows that one
cannot use the Ptolemaic criteria for the evaluation of achievements made on
the basis of the Newtonian approach. Progress in the development of the
fundamental theory is connected mainly with the reduction of the number of
the fundamental statements.

Existence of the hidden mistake in the fundamental statements of a theory, a
use of the Ptolemaic approach to the development of the theory and Ptolemaic
criteria in the evaluation of the extent of this development lead to the
crisis, i.e. to the surcease of the theory development. Well known example
of such a situation is the interplay between the Ptolemaic doctrine and the
Copernicus one. The origin of the crisis was the mistake in the
determination of the mutual motion of the Sun and the Earth. But the mistake
in itself could not generate the crisis. The crisis appeared, because
advocates of Ptolemeus and those of Copernicus used different methods of
evaluation of the theory quality. As a result the Ptolemaic advocates did
not accept the correction suggested by Copernicus.

Another example: the crisis in the foundation of geometry, when in the
beginning of XIXth century the non-Euclidean geometry was suggested by N.
Lobachewsky and J. Bolyai. Mathematicians of that time did not want to
consider the non-Euclidean geometry as a possible geometry of the real
world. It is not clear exactly, what was the origin of the rejection of the
non-Euclidean geometry. I do not know investigations concerning this
problem. My opinion in this question is as follows.

One can always introduce the Cartesian coordinate system in the Euclidean
geometry. In the Riemannian geometry it was possible not always. This
circumstance was the main argument against application of the Riemannian
geometry to the real world. Most of mathematicians of that time considered
the Cartesian coordinate system to be an attribute of any geometry, and
inapplicable properties of the coordinate system were considered as
inapplicable properties of the geometry in question. The crisis was not so
distinct, as in the case Ptolemeus - Copernicus, but it was a crisis,
because a rejection of a true theory is one of symptoms of the crisis. In
the beginning of the XXth century the non-Euclidean (Riemannian) geometry
was accepted as a possible space-time geometry, and it was clear that the
rejection of the non-Euclidean geometry was erroneous. In general, the
crisis is a social phenomenon. A mistake in the fundamental statements of a
fundamental theory is not yet a crisis. The crises appears, when the most
researchers use an erroneous theory, rejecting the true one.

Finally, the present crisis in the contemporary geometry is connected with
the incorrect generalization of properties of the Euclidean geometry. The
contemporary geometers consider that the straight is the one-dimensional
line in all geometries and, in particular, in the real space-time geometry.
The Riemannian space-time geometry is based on this erroneous statement, and
insistence on this statement is an origin of the crisis in the contemporary
geometry. The incomprehensible rejection by editors of mathematical
journals, cited above, is a manifestation of this crisis. Unfortunately, a
logical motivation does not help in the overcoming of the crisis situation,
because by definition the crisis is such a situation, where the usual logic
cannot overcome the preconception.

The main preconception is the fundamental supposition that the straight is a
one-dimensional line in all possible geometries. This preconception is
followed by the supposition, that the one-dimensional curve is a primary
geometrical object, which is to be used essentially at the construction of
any geometry. These statements seem to be very reasonable. As well as other
geometers, I was holding this viewpoint for a long time. Why was I forced to
reject these reasonable ideas? Maybe, a history how I succeeded to overcome
this preconception helps to other researchers in overcoming of this
preconception. Here I present a review of ideas of the tubular geometry
(T-geometry), where the supposition on one-dimensionality of the straight
and that on the leading role of the curve at the geometry construction are
violated. Of necessity my presentation will be a review of my own papers,
because I do not know other researchers, participating in the creation of
the T-geometry. Besides, the autobiographical style of the review allows one
to present also motives of my investigations, which cannot be found in the
published papers, where the presentation is produced formally. The knowledge
of motives may appear to be useful in the overcoming of the crisis situation.

\section{The Riemannian geometry in terms of the world function}

My first paper was written in the end of 1958, submitted in 1959 and
published in 1962 \cite{R62}. Writing the paper, I was the fifth-year
student of the Moscow State University. In that time I was interested in the
Riemannian geometry and in the Einstein theory of gravitation.
Conventionally the Riemannian geometry is described in terms of the
infinitesimal space-time interval $dS$ and metric tensor $g_{ik}$

\begin{equation}
dS^{2}=g_{ik}dx^{i}dx^{k}  \label{a2.1}
\end{equation}%
and the geodesic is described in terms of ordinary differential equations.

If one introduces the finite interval 
\begin{equation}
S\left( x,x^{\prime }\right) =\int_{\mathcal{L}_{xx^{\prime }}}\sqrt{%
g_{ik}dx^{i}dx^{k}},\qquad \sigma \left( x,x^{\prime }\right) =\frac{1}{2}%
S^{2}\left( x,x^{\prime }\right)  \label{a2.2}
\end{equation}%
where integration is produced along the geodesic $\mathcal{L}_{xx^{\prime }}$%
, connecting points $x$ and $x^{\prime }$, the finite interval $S\left(
x,x^{\prime }\right) $ satisfies the Jacobi-Hamilton equation 
\begin{equation}
\frac{\partial S}{\partial x^{i}}g^{ik}\frac{\partial S}{\partial x^{k}}=1
\label{a2.3}
\end{equation}%
In terms of the world function $\sigma \left( x,x^{\prime }\right) =\frac{1}{%
2}S^{2}\left( x,x^{\prime }\right) $ this equation takes the form 
\begin{equation}
\sigma _{i}g^{ik}\left( x\right) \sigma _{k}=2\sigma ,\qquad \sigma
_{i}\equiv \frac{\partial \sigma }{\partial x^{i}}  \label{a2.4}
\end{equation}%
In terms of $\sigma $ the geodesic $x^{k}=x^{k}\left( \tau \right) $ is
described by the algebraic equations 
\begin{equation}
\sigma _{k^{\prime }}\left( x,x^{\prime }\right) =b_{k^{\prime }}\tau
,\qquad \sigma _{k^{\prime }}\left( x,x^{\prime }\right) \equiv \frac{%
\partial \sigma }{\partial x^{\prime k^{\prime }}},\qquad b_{k^{\prime }}=%
\text{const},\qquad k^{\prime }=0,1,...n  \label{a2.5}
\end{equation}%
where $b_{k^{\prime }}$ is a vector at the point $x^{\prime }$, and the
geodesic $x^{k}=x^{k}\left( \tau \right) $ is determined implicitly by
algebraic equations (\ref{a2.5}). If (\ref{a2.2}) is considered to be the
action functional, describing the motion of a free relativistic particle,
the equation (\ref{a2.3}) is the Jacobi-Hamilton equation for this action,
and (\ref{a2.5}) describes an obtaining of the particle world line from the
solution of the Jacobi-Hamilton equation. The equations (\ref{a2.5})
determine the geodesic, passing through the point $x^{\prime }$ in the
direction of the vector $b_{k^{\prime }}$.

It follows from (\ref{a2.5}) that the description in terms of the world
function $\sigma $ is more informative, than the description in terms of the
metric tensor $g_{ik}$, because transition from the metric tensor $g_{ik}$
to the world function $\sigma $ leads to a replacement of the differential
equations by the algebraic ones.

In the paper \cite{R62} (English version of this presentation in \cite{R64})
I succeeded to construct the mathematical formalism, which admits one to
describe the $n$-dimensional Riemannian geometry in terms of the world
function $\sigma $, given on the $n$-dimensional manifold. It is very
important that a system of differential equations for the world function had
been obtained. This system of differential equations contains only the $%
\sigma $ function and its derivatives 
\begin{eqnarray}
\sigma _{i}\sigma ^{i,k^{\prime }}\sigma _{k^{\prime }} &=&2\sigma ,\qquad
\sigma (x,x^{\prime })=\sigma (x^{\prime },x),\qquad \sigma (x,x)=0
\label{a2.9} \\
G_{ik||l} &=&0,\qquad i,k,l=1,2,...,n  \label{a2.10}
\end{eqnarray}%
where $\sigma ^{i,k^{\prime }}$ is determined by 
\begin{equation}
\sigma ^{i,k^{\prime }}\sigma _{l,k^{\prime }}=\delta _{l}^{i},\qquad
i,l=1,2,...,n  \label{a2.11}
\end{equation}%
\begin{equation}
\sigma _{l,k^{\prime }}\equiv \frac{\partial ^{2}\sigma (x,x^{\prime })}{%
\partial x^{l}\partial x^{\prime k^{\prime }}},\qquad k^{\prime },l=1,2,...,n
\label{a2.12}
\end{equation}%
The symbol $(.)_{||l}$ denotes the tangent derivative with respect to $x^{l}$%
, i.e. the covariant derivative with respect to $x^{l}$ with the Christoffel
symbol 
\begin{equation}
\Gamma _{kl}^{i}(x,x^{\prime })=\sigma ^{i,j^{\prime }}\sigma _{klj^{\prime
}},\qquad \sigma _{klj^{\prime }}\equiv \frac{\partial ^{3}\sigma
(x,x^{\prime })}{\partial x^{k}\partial x^{l}\partial x^{\prime j^{\prime }}}%
,\qquad i,k,l,j^{\prime }=1,2,...,n  \label{a2.13}
\end{equation}%
and $G_{ik}$ is defined by the relation 
\begin{equation}
G_{ik}=G_{ik}(x,x^{\prime })\equiv \sigma _{i||k}=\frac{\partial \sigma _{i}%
}{\partial x^{k}}-\Gamma _{ik}^{I}\sigma _{l},\qquad i,k=1,2,...,n
\label{a2.14}
\end{equation}

The conditions (\ref{a2.9}), (\ref{a2.10}) are the necessary and sufficient
conditions of that the world function $\sigma $ describes the Riemannian
geometry. The system (\ref{a2.9}), (\ref{a2.10}) allows one to reject
definition (\ref{a2.2}) of the world function $\sigma $ which defines $%
\sigma $ via the metric tensor $g_{ik}$ and the concept of the curve
(geodesic $\mathcal{L}_{xx^{\prime }}$). Now the world function $\sigma $
may be defined as a symmetric function $\sigma \left( x,x^{\prime }\right)
=\sigma \left( x^{\prime },x\right) $, satisfying the system of differential
equations (\ref{a2.9}), (\ref{a2.10}).

Such a definition of the world function put the following question. Let us
imagine, that the function $\sigma $ does not satisfy the system (\ref{a2.9}%
), (\ref{a2.10}) of differential equations. What does the world function
formalism describe in this case: a non-Riemannian geometry, or no geometry?
It was a very difficult question, and I could not answer it for a long time.
In this case the world function cannot be defined by the relation (\ref{a2.2}%
).

In the Riemannian geometry the segment $\mathcal{L}_{[xx^{\prime }]}$ of a
geodesic $\mathcal{L}_{xx^{\prime }}$ between the points $x$ and $x^{\prime
} $ can be defined directly via the finite interval $S\left( x,x^{\prime
}\right) =\sqrt{2\sigma \left( x,x^{\prime }\right) }$ by means of the
relation 
\begin{equation}
\mathcal{L}_{[xx^{\prime }]}=\left\{ z|S\left( x,z\right) +S\left(
z,x^{\prime }\right) =S\left( x,x^{\prime }\right) \right\}  \label{a2.15}
\end{equation}%
The relation (\ref{a2.15}) defines the geodesic as a set of points without a
reference to the concept of a curve. The definition (\ref{a2.15}) has
another advantage. It does not contain a reference to a coordinate system.
Although the relation (\ref{a2.15}) contains a labelling of points by the
coordinates, but it can be easily written in the form 
\begin{equation}
\mathcal{L}_{[PQ]}=\left\{ R|S\left( P,R\right) +S\left( R,Q\right) =S\left(
P,Q\right) \right\} ,\qquad P,Q,R\in \Omega  \label{a2.16}
\end{equation}%
where $\Omega $ is a set of points, where the geometry is given.

However, the relation (\ref{a2.15}) describes a multi-dimensional surface,
which, in general, does not degenerate into a one-dimensional line in the
case, when the world function does not satisfy constraints (\ref{a2.9}), (%
\ref{a2.10}). This circumstance hinders from an introduction of the
non-Riemannian geometry, because I could not imagine a geometry with
surfaces (tubes) instead of straight lines. I think the idea, that the
straight is a one-dimensional line in all geometries, hinders other
researchers from construction of the non-Riemannian geometry. Such an idea
is an associative delusion, which reminds the delusion of the ancient
Egyptians, who believed that all rivers flow towards the North.

When my paper \cite{R62} was in print (1959 - 1962), the excellent book by
J.L. Synge \cite{S60} appeared in libraries of USSR. In this book the
Riemannian geometry and the general relativity were presented in terms of
the world function. Writing the paper \cite{R62}, I knew nothing on the
papers of Synge, devoted to the description of the Riemannian geometry in
terms of the world function. I constructed the world function formalism from
the outset, and my formalism distinguished essentially from the formalism of
Synge. The Synge's formalism was oriented to the expansion of the world
function $\sigma \left( x,x^{\prime }\right) $ over powers of the $%
x-x^{\prime }$, whereas my formalism was oriented to such an expansion in
the minimal degree. My formalism used the world function as a means of
description, which may be an alternative to the coordinate description and
to the use of the metric tensor. In particular, the book of Synge does not
contain equations (\ref{a2.9}), (\ref{a2.10}), which were considered by me
as the principal result of my paper. Comparison of my paper and the Synge's
book impressed me (I was a student, whereas J.L. Synge was the famous
researcher) and changed my self-concept. Starting from the same position as
J.L. Synge, I have obtained in half-year the more developed results. It
meant that I possessed a nonstandard way of thinking. It meant that I should
avoid a detailed study of the literature before the investigation, because
together with results of investigation I risked to perceive the
preconceptions. In the case of the crisis situation in physics and geometry
an elimination of preconceptions is much more important, than obtaining of
the positive results of investigations. Results of correct investigations
could be repeated easily, but if the result of other researchers has been
obtained on the basis of some preconception, they influence upon my way of
thinking, and the discovery of the preconception appears to be a more
difficult problem. Of course, such an approach is justified only in the
crisis situation, when ,on the one hand, one is sure that the hidden
preconceptions and mistakes take place, and on the other hand, one is sure
in his own capacities.

\section{Further development of T-geometry}

I have returned to problems of geometry only in 1989. It was clear for me,
that the description in terms of the world function was a very effective
method of the geometry description. I wanted to refine maximally the
geometry description from influence of the other methods of description. In
particular, I wanted to obtain a coordinateless description, as it has been
made by Euclid many years ago. In general, I believed that the geometry was
a totality of geometrical objects and relations between them, but not a
totality of coordinates, metric tensor, curvature tensor, covariant
derivatives, Christoffel symbols, etc. I understood, that the coordinates
were not a necessary attribute of a geometry, but coordinates were a
necessary attribute of dynamics. At least, I did not know examples of a
coordinateless description of dynamics (but the example of the
coordinateless description of the Euclidean geometry was known). Describing
dynamics of particles in the space-time, we cannot avoid a use of
coordinates in dynamics and in the space-time geometry. Nevertheless, the
geometry should be constructed in the coordinateless form, to separate the
properties of the geometry from the properties of the description.

I tried to present the geometry in the form, which be maximally free of the
means of description. To make this, one should understand what is the
geometry in itself and what is the means of the description. The world
function $\sigma $ as a means of the description is better, than the metric
tensor $g_{ik}$, because $\sigma $ is a scalar, i.e. the quantity invariant
with respect to the coordinate transformation, whereas the metric tensor
does not possess this property. The Euclidean straight $\mathcal{T}%
_{P_{0}P_{1}}$, passing through the points $P_{1}$, $P_{2}$ can be expressed
via the world function by means of the relation 
\begin{equation}
\mathcal{T}_{P_{0}P_{1}}=\left\{ R|F_{2}\left( P_{0},P_{2},R\right)
=0\right\}  \label{a3.4}
\end{equation}
where 
\begin{equation}
F_{n}\left( \mathcal{P}^{n}\right) =\det ||(\mathbf{P}_{0}\mathbf{P}_{i}.%
\mathbf{P}_{0}\mathbf{P}_{k})||,\qquad i,k=1,2,...n,\qquad \mathcal{P}%
^{n}\equiv \left\{ P_{0},P_{1},...P_{n}\right\}  \label{a3.5}
\end{equation}
is the $n$th order Gram determinant, constructed of scalar products 
\begin{equation}
(\mathbf{P}_{0}\mathbf{P}_{i}.\mathbf{P}_{0}\mathbf{P}_{k})=\sigma \left(
P_{0},P_{i}\right) +\sigma \left( P_{0},P_{k}\right) -\sigma \left(
P_{i},P_{k}\right)  \label{a3.6}
\end{equation}
In the proper Euclidean geometry the relation (\ref{a3.6}) is a corollary of
evident relations 
\begin{equation}
\left\vert \mathbf{P}_{0}\mathbf{P}_{i}\right\vert ^{2}=2\sigma \left(
P_{0},P_{i}\right)  \label{a3.6a}
\end{equation}
\begin{equation}
\left\vert \mathbf{P}_{i}\mathbf{P}_{k}\right\vert ^{2}=\left\vert \mathbf{P}%
_{0}\mathbf{P}_{k}-\mathbf{P}_{0}\mathbf{P}_{i}\right\vert ^{2}=\left\vert 
\mathbf{P}_{0}\mathbf{P}_{i}\right\vert ^{2}+\left\vert \mathbf{P}_{0}%
\mathbf{P}_{k}\right\vert ^{2}-2(\mathbf{P}_{0}\mathbf{P}_{i}.\mathbf{P}_{0}%
\mathbf{P}_{k})  \label{a3.6b}
\end{equation}

In the proper Euclidean space the condition of the Gram determinant
disappearance

\begin{equation}
F_{n}\left( \mathcal{P}^{n}\right) =0  \label{a3.7}
\end{equation}%
is the necessary and sufficient condition of the linear dependence of $n$
vectors $\mathbf{P}_{0}\mathbf{P}_{1}$, $\mathbf{P}_{0}\mathbf{P}_{2}$,...$%
\mathbf{P}_{0}\mathbf{P}_{n}$. If the world function of the Euclidean space
is slightly changed, the set of points (\ref{a3.4}) ceases to be the
one-dimensional line and turns into a surface, whose dimension is more than $%
1$. Besides, in the Minkowski space-time geometry (four-dimensional
pseudo-Euclidean geometry of index 1) the timelike straight $\mathcal{T}%
_{P_{0}P_{1}}$, ($\sigma \left( P_{0},P_{1}\right) >0$), i.e. the set of
points (\ref{a3.4}) is a straight line, whereas the spacelike straight $%
\mathcal{T}_{P_{0}P_{1}}$, ($\sigma \left( P_{0},P_{1}\right) <0$, the set
of points (\ref{a3.4})) is a three-dimensional surface. Such a situation
appears, because in the Minkowski geometry the condition (\ref{a3.7}), (\ref%
{a3.5}) does not coincide with the condition of the linear dependence of $n$
vectors $\mathbf{P}_{0}\mathbf{P}_{1}$, $\mathbf{P}_{0}\mathbf{P}_{2}$,...$%
\mathbf{P}_{0}\mathbf{P}_{n}$ in its conventional form. Conventionally the
linear dependence of vectors in the Minkowski geometry is determined on the
basis of the linear vector space, which may be defined in the Minkowski
space.

It seems that the geometry construction in the Minkowski space depends on
the way how the linear dependence of vectors is introduced. In other words,
the geometry depends on the axioms, which are used at its construction. It
is true, provided the geometry is considered simply as a logical
construction. But here we deal with the physical geometry, which pretends to
be a geometry of the real space-time, which has not to depend on our choice
of basic axioms. The axioms, if we use them at construction of the Minkowski
geometry, must be given objectively, i.e. in such a way, that they agree
with experimental data. Practically it means, that we must construct a
maximally broad class of geometries and thereafter choose the geometry,
which satisfies the experimental data. The condition (\ref{a3.7}), (\ref%
{a3.5}) provides construction of maximally general geometry, because it does
not depend on any constraints on geometry.

The idea of the tubular character of straights was less exotic for
physicists, than for mathematicians. Dealing with the microcosm, physicists
invent essentially more exotic objects: such as strings, branes and
stochastic geometries. However, these objects are external with respect to
the geometry, because they do not demand the geometry reconstructions. The
geometry existed in itself, and the strings and other exotic objects existed
in themselves without interaction with the geometry. The situation with the
straight was other. One can imagine the thin hallow tube instead of the
straight. In the large scale the tube thickness was imperceptible, and no
new effects connected with the tube thickness appeared. The straight was the
primary object (but not the external one), i.e. the construction of the
conventional Riemannian geometry was based on the properties of the
straight. Any change of the properties of the straight generated a
reconstruction of the whole geometry. My brain slept until it became clear
for me, that the tubular character of straight allows one to construct such
a space-time geometry, where the free particle motion were primordially
stochastic. Such a geometry could replace the principles of quantum
mechanics. 

I admit that other researchers came to the idea of the tubular geometry. But
what to do with such a geometry? What profit can one obtain from such a
geometry? The pure geometers could construct the tubular geometry, but they
did not know what to do with it. In distinction to them I was a physicist,
and I knew how to use such a geometry. Several years before the T-geometry
construction I succeeded to construct the \textquotedblright relativistic
statistics\textquotedblright , i.e. the conception, where the
nonrelativistic quantum mechanics was obtained as a result of statistical
description of stochastically moving particles. This conception had only one
defect. It was unclear why the free particles moved stochastically. The
T-geometry explained freely the origin of the stochasticity by the
space-time properties. All further questions disappeared at once. Thus,
appearance of the T-geometry is explained by its necessity. But only the
physicist, who succeeded to construct the \textquotedblright relativistic
statistics\textquotedblright , perceived this necessity.

Let us list arguments in favor of the definition (\ref{a3.7}), (\ref{a3.5})
of the linear dependence:

\begin{enumerate}
\item The definition of the linear dependence (\ref{a3.7}), (\ref{a3.5}) may
be used on any set $\Omega $ of points, where the world function $\sigma $
is given. This set $\Omega $ may be continuous or discrete. It may be
inhomogeneous and non-isotropic

\item The definition (\ref{a3.7}), (\ref{a3.5}) does not refer to a
coordinate system, or some other additional structure, which should be
introduced.

\item The definition (\ref{a3.7}), (\ref{a3.5}) agrees with the experimental
data: the timelike one-dimensional straight describes a motion of a free
particle. The spacelike straight does not associate with the particle
motion, or associates with motion of the hypothetical tachyon. It has not
been discovered, and nobody knows, whether tachyon is described by a
spacelike one-dimensional straight, or by a three-dimensional surface.
\end{enumerate}

The conventional definition, based on the introduction of the linear vector
space, has the advantage, that it was used before, and it is customary for
all researchers. Unfortunately, the conventional definition contains a
reference to the additional structure (linear vector space), which can be
introduced only for uniform and continuous geometries. As a result the class
of geometries, where the conventional definition of the linear dependence
can be used, appears to be rather narrow, and this class does not contain
the real space-time geometry.

At first, my preconception against the definition (\ref{a3.7}), (\ref{a3.5})
were very strong. On the other hand, I wanted to describe the Riemannian
geometry and its possible generalization in the coordinateless form. It was
possible only in terms of the world function. My solution of this question
was half-and-half. I assumed that the relation (\ref{a3.4}) describes the
geometrical object $\mathcal{T}_{P_{0}P_{1}}$, consisting of one-dimensional
lines (the line tube), and for determination of one line one needs to impose
additional constraints. The additional constraints, which should be imposed
on the set $\mathcal{T}_{P_{0}P_{1}}$ described properties of the Riemannian
geodesic.

In the conventional approach a geodesic in a $n$-dimensional Riemannian
space is considered to be a \textit{special type of a curve} having the
following extremal properties

(i) \textit{Extremality.} The distance $(2\sigma )^{1/2}$ between two points
measured along a geodesic is the shortest (extremal) as compared with the
distance measured along other curves.

(ii) \textit{Definiteness. }Any two points of the geodesic determine
unambiguously the geodesic passing through these points.

(iii) \textit{Minimality of section} (one-dimensionality). Any section of a
geodesic consists of one point.

Another approach appears, when the geodesic is considered to be a \textit{%
special type of a surface} (or a line tube) that degenerates into a line.
Then properties (ii) and (iii) appear to be fulfilled; however, the property
(i) is not defined because the concept of a line (or curve) is not defined.
It appeared, that the extremal property (i) is equivalent to properties (ii)
and (iii), taken together with (\ref{a3.4}). Let us try to define a geodesic
tube having the property of definiteness (ii) and that of minimality of
section (iii) at the same time. In the paper \cite{R90} the property (ii)
and (iii) were written in terms of the world function in the form

(ii) definiteness

\begin{equation}
\forall Q_{0},Q_{1}\in \mathcal{T}_{P_{0}P_{1}}\wedge \sigma \left(
Q_{0},Q_{1}\right) \neq 0,\qquad \mathcal{T}_{Q_{0}Q_{1}}\subset \mathcal{T}%
_{P_{0}P_{1}}  \label{a3.8}
\end{equation}%
(iii) minimality of the section 
\begin{equation}
P\in \mathcal{T}_{P_{0}P_{1}},\qquad \mathcal{S}_{1,P}\left( \mathcal{T}%
_{P_{0}P_{1}}\right) \equiv \left\{ R|\dbigwedge\limits_{k=0,1}\sigma \left(
P_{k},P\right) =\sigma \left( P_{k},R\right) \right\} =\left\{ P\right\} 
\label{a3.9}
\end{equation}%
The properties (ii) and (iii) together with (\ref{a3.4}) were used instead
of the property (i), which is not defined for tubes determined by the
relation (\ref{a3.4}). Combination of relations (\ref{a3.8}), (\ref{a3.9})
with the tube equation (\ref{a3.4}) and subsequent investigation in \cite%
{R90} allow one to study the properties of the curve, described in terms of
the world function. It was rather strange and unexpected that the tube (\ref%
{a3.4}) and the high order tubes are primary geometrical objects, but not a
unification of many one-dimensional lines. It was strange and unexpected
that the one-dimensional line was a derivative geometrical object, defined
as an intersection of several surfaces, which were considered to be the
primary objects (tubes). It was strange and not customary, because at the
conventional approach to the Riemannian geometry any surface may be
considered as consisting of many primary geometrical objects --
one-dimensional lines.

Thus at the conventional approach the curve is a primary object, which lies
in the foundation of the Riemannian geometry, and the world function is
defined by means of (\ref{a2.2}) via the concept of the curve. Such a
definition may lead to an inconsistency. Let me explain this.

Euclid defined a curve line as a limit of the broken line $\mathcal{T}_{%
\mathrm{br}}$, consisting of many rectilinear segments $\mathcal{T}_{\left[
P_{k}P_{k+1}\right] }$ 
\begin{equation}
\mathcal{T}_{\mathrm{br}}=\dbigcup\limits_{k}\mathcal{T}_{\left[ P_{k}P_{k+1}%
\right] }  \label{a3.1}
\end{equation}%
The curve appeared, when the length of any segment $\mathcal{T}_{\left[
P_{k}P_{k+1}\right] }$ tended to zero. But even if the limit is not taken,
the broken line (\ref{a3.1}) is a kind of a curve, which is not smooth, in
general. Thus, Euclid defined the curve via the straight. The straight was
the primary geometrical object, whereas the curve was the secondary
(derivative) geometrical object. It was the metrical definition of the curve.

In the Riemannian geometry one uses the topological (non-metrical)
definition of the curve. The curve $\mathcal{L}$ is defined as a continuous
mapping 
\begin{equation}
\mathcal{L}:\qquad \left[ 0,1\right] \rightarrow \Omega ,\qquad \left[ 0,1%
\right] \subset \mathbb{R}  \label{a3.2}
\end{equation}
The straight line (geodesic) $\mathcal{T}_{P_{0}P_{1}}$ between the points $%
P_{0}$ and $P_{1}$ is defined as a special case of the curve $\mathcal{L}$,
when the length of the curve between the points $P_{0}$ and $P_{1}$ is
minimal. Thus, in the Riemannian geometry the curve $\mathcal{L}$ is the
primary geometrical object, whereas the straight is a derivative (secondary)
geometrical object. The world function, defined by the relation (\ref{a2.2})
is also a derivative quantity. The definition of the curve (\ref{a3.2})
needs an introduction of additional concepts such as topological space,
continuity, etc.

Compatibility of two different definitions (\ref{a3.1}) and (\ref{a3.2}) is
not evident. If we demand the coincidence of the definitions (\ref{a3.1})
and (\ref{a3.2}), we may overdetermine the problem of the geometry
construction. As a result we may discriminate some geometries, which are
interesting as the space-time geometries. If it appears that the geometry
can be constructed on the basis of only the world function and of the
metrical definition of the curve (\ref{a3.1}), the topological definition (%
\ref{a3.2}) becomes to be redundant. Using them as primary concepts at the
construction of geometry, one may lead to inconsistencies.

In some cases the definition (\ref{a2.2}) leads to absurd results. Let $%
\mathcal{L}_{2}$ be the two-dimensional proper Euclidean plane, and $K$ be
the Cartesian coordinate system on $\mathcal{L}_{2}$. If the world function
is considered to be a primary quantity, it is described by the formula 
\begin{equation}
\sigma \left( \mathbf{r},\mathbf{r}^{\prime }\right) =\frac{1}{2}\left\vert 
\mathbf{r}-\mathbf{r}^{\prime }\right\vert ^{2},\qquad \mathbf{r=}\left\{
x,y\right\} ,\qquad \mathbf{r}^{\prime }\mathbf{=}\left\{ x^{\prime
},y^{\prime }\right\}   \label{a3.10}
\end{equation}%
If the proper Euclidean geometry on $\mathcal{L}_{2}$ is considered as a
special case of the Riemannian geometry, the world function is determined by
the relation (\ref{a2.2}). The world function has the form (\ref{a3.10}),
provided we consider the whole plane $\mathcal{L}_{2}$. Now let us cut the
plane $\mathcal{L}_{2}$ into two parts $\mathcal{L}_{2c}=\left\{ \mathbf{r|r}%
^{2}\leq 1\right\} $ and $\mathcal{L}_{2o}=\left\{ \mathbf{r|r}%
^{2}>1\right\} $ and define the world function in $\mathcal{L}_{2c}$ and in $%
\mathcal{L}_{2o}$, using (\ref{a2.2}). In $\mathcal{L}_{2c}$ we obtain the
result (\ref{a3.10}), whereas in $\mathcal{L}_{2o}$ we obtain the world
function $\sigma \left( P,Q\right) $, distinguishing from (\ref{a3.10}), if
the straight line segment $\mathcal{T}_{[PQ]},$ connecting the points $P$
and $Q$, crosses $\mathcal{L}_{2c}$. We obtain the absurd result that the
region $\mathcal{L}_{2o}$ cannot be embedded into $\mathcal{L}_{2}$
isometrically. We obtain the convexity problem, when in the region $\mathcal{%
L}_{\mathrm{cut}}$, which is cut out of the space $\mathcal{L}$, the
Riemannian geometry is violated, if the region $\mathcal{L}_{\mathrm{cut}}$
is nonconvex. 

If the world function is a primary quantity, the geometry of $\mathcal{L}_{%
\mathrm{cut}}$ is conserved independently of convexity or nonconvexity of $%
\mathcal{L}_{\mathrm{cut}}$. It means that the \textit{convexity problem is
an artificial problem}, connected with application of the topological
concept of a curve (\ref{a3.2}) in the construction of the Riemannian
geometry. As far as the Riemannian geometry may be constructed without a
reference to the topological concept of a curve, and without a reference to
topological concepts, the topological concepts appear to be only derivative
concepts.

The next step in the development of the T-geometry was made, when one
discovered that the T-geometry allowed one to choose the space-time geometry
in such a way, that the motion of free particles appeared to be primordially
stochastic, and statistical description of this stochastic motion was
equivalent to the quantum description in terms of the Schr\"{o}dinger
equation \cite{R91}. According to the Newtonian criterion it was a great
progress of the fundamental theory, because the quantum mechanics principles
appeared to be a corollaries of the space-time geometry and of the
statistical description. They ceased to be additional hypotheses, invented
for explanation of quantum phenomena.

According to the Ptolemaic criterion this result was not a progress, until
it will be shown that it can explain any new physical phenomena. It
presented only another version of the quantum phenomena explanation, and
this new version was not better, the old version of the quantum mechanics,
because it did not explain new effects. From the viewpoint of the Ptolemaic
criterion the reduction of the number of fundamental principles was of no
importance.

As to me, I stayed on the viewpoint of Newton and the successful application
of T-geometry to the space-time confirmed validity of my approach to
geometry. The obtained world function for the real space-time has the form 
\begin{equation}
\sigma =\sigma _{\mathrm{M}}+\left\{ 
\begin{array}{c}
\frac{\hbar }{2bc},\text{ if }\sigma _{\mathrm{M}}>\sigma _{0} \\ 
0,\text{ if }\sigma _{\mathrm{M}}<0%
\end{array}
\right.  \label{a3.11}
\end{equation}
where $\sigma _{\mathrm{M}}$ is the world function of the Minkowski space, $%
\hbar $ is the quantum constant, $c$ is the speed of the light, $b\lesssim
10^{-17}$g/cm is some universal constant an the quantity $\sigma _{0}\approx 
\frac{\hbar }{2bc}$. The world function $\sigma $ is not determined for $%
0<\sigma _{\mathrm{M}}<\sigma _{0}$. In other words, the value of the world
function $\sigma $ for $0<\sigma _{\mathrm{M}}<\sigma _{0}$ is inessential
for description of quantum phenomena. On the other hand, such quantities as
the metric tensor and the curvature tensor are defined as the expansion
coefficients of the $\sigma \left( x,x^{\prime }\right) $ expansion over
powers of $x-x^{\prime }$, i.e. in the vicinity of $\sigma _{\mathrm{M}}=0$.
As far as the value of $\sigma $ at $x-x^{\prime }=0$ (or $\sigma _{\mathrm{M%
}}=0$) are inessential for description of quantum phenomena the values of
the metric tensor and of the curvature tensor appears to be inessential for
description of quantum phenomena. The conventional local description of the
Riemannian geometry, as well as the world function technique on the
manifold, developed in \cite{R62} do not work in the nondegenerate geometry
with tubes instead of straights. It became clear after my attempt of
application of the methods of differential geometry to the nondegenerate
geometry \cite{R92}.

In the paper \cite{R90} a coordinateless description of the Riemannian
geometry was considered, where the world function was a tool of this
description. Such a description generates some generalization of the
Riemannian geometry. But the term T-geometry was not mentioned in the paper 
\cite{R90}, and the role of the deformations of the Euclidean geometrical
objects was not clear for me, although I used deformations of the Euclidean
straights and planes practically.

\section{T-geometry as a generalization of metric \\ geometry}

Ten years ago I returned to problems of the T-geometry, but then the
T-geometry was considered as the generalization of the metric geometry,
whereas in the papers \cite{R90,R91,R92} it was considered as a
generalization of the Riemannian geometry. What is the difference between
the two approaches?

The metric geometry is given on the metric space, which is described by the
metric $\rho $. The metric space $M=\{\rho ,\Omega \}$ is a set $\Omega $ of
points $P\in \Omega $ with the metric $\rho $ given on $\Omega \times \Omega 
$ 
\begin{equation}
\rho :\quad \Omega \times \Omega \rightarrow D_{+}\subset \mathbb{R}
\label{a1.1}
\end{equation}%
\begin{equation}
\rho (P,P)=0,\qquad \rho (P,Q)=\rho (Q,P),\qquad \forall P,Q\in \Omega
\label{a1.2}
\end{equation}%
\begin{equation}
D_{+}=[0,\infty ),\qquad \rho (P,Q)=0,\quad \text{if and only if}\mathrm{%
\;\;\;\;}P=Q,\qquad \forall P,Q\in \Omega  \label{a1.3}
\end{equation}%
\begin{equation}
\rho (P,Q)+\rho (Q,R)\geq \rho (P,R),\qquad \forall P,Q,R\in \Omega
\label{a1.4}
\end{equation}%
The T-geometry is defined on $\sigma $-space $V=\left\{ \sigma ,\Omega
\right\} $ by means of relations 
\begin{equation}
\sigma :\quad \Omega \times \Omega \rightarrow \mathbb{R}  \label{a3.14}
\end{equation}%
\begin{equation}
\sigma (P,P)=0,\qquad \sigma (P,Q)=\sigma (Q,P),\qquad \forall P,Q\in \Omega
\label{a3.15}
\end{equation}

It seems that the metric geometry is a special case of the T-geometry, when $%
\sigma =\frac{1}{2}\rho ^{2}$, and the metric is restricted by constraints (%
\ref{a1.3}), (\ref{a1.4}). But it is not so.

What is actually happens is that the T-geometry is equipped by the
deformation principle, whereas the metric geometry is free of this
equipment. The deformation principle means as follows. Any geometrical
object $\mathcal{O}$ and any relation $\mathcal{R}$ between the geometrical
objects of the proper Euclidean geometry can be expressed in terms of the
Euclidean world function $\sigma _{\mathrm{E}}$. To obtain expressions for $%
\mathcal{O}$ and $\mathcal{R}$ in the arbitrary T-geometry, it is sufficient
to replace $\sigma _{\mathrm{E}}\rightarrow \sigma $ in the Euclidean
expressions for $\mathcal{O}$ and $\mathcal{R}$. Of course, one can easily
equip the metric space by the deformation principle, but it had not been
made somehow.

I think that this fact takes place, because of the underestimate of the
deformation principle role in the construction of a geometry. In the
Riemannian geometry the deformation of the Euclidean space is used, when the
Euclidean infinitesimal interval $dS_{\mathrm{E}}^{2}=g_{\mathrm{E}%
ik}dx^{i}dx^{k}$ is replaced by the Riemannian one $dS_{\mathrm{R}}^{2}=g_{%
\mathrm{R}ik}dx^{i}dx^{k}$. It is this replacement, that determines the
essence of the obtained Riemannian geometry. The features and
characteristics of this change are to be investigated, if we want to study
the Riemannian geometry. Unfortunately, comparison of $dS_{\mathrm{E}}^{2}$
and $dS_{\mathrm{R}}^{2}$ is very complicated technically, because one needs
to compare metric tensors and coordinate systems simultaneously. In most
textbooks on the Riemannian geometry most attention has been concentrated on
technical details (mappings, coordinate transformations, covariant
differentiation and other means of description), whereas the conceptual
procedure of replacement $dS_{\mathrm{E}}^{2}$ by $dS_{\mathrm{R}}^{2}$
appears outside of the field of vision. As a result most of geometers
believe, that the Riemannian geometry is a totality of coordinates, metric
tensors, covariant derivatives and other means of the Riemannian geometry
description. They appear to be discouraged, when they discover that the
geometry and the means of the geometry description are different things.

Only by this circumstance one can explain the fact that, constructing the
metric geometry, one considers only the metric, but forgets about the
deformation principle. Instead of the deformation principle one uses
additional constraints (\ref{a1.3}), (\ref{a1.4}). Using the metric
definition of the curve (\ref{a3.1}) and constraints (\ref{a1.3}), (\ref%
{a1.4}), providing one-dimensional character of segments $\mathcal{T}%
_{[P_{i}P_{k}]}$, one succeeded to introduce the concept of an
one-dimensional curve in the metric geometry. The shortest line between two
points may be considered as the straight line (the shortest) in the metric
geometry. But already an introduction of two-dimensional plane in the metric
geometry meets difficulties. As a result, the metric geometry appears to be
a very poor geometry. K. Menger \cite{M28} and L.M. Blumental \cite{B53}
tried to construct the distance geometry, which was the metric geometry
without the constraints (\ref{a1.3}), (\ref{a1.4}). They did not use the
deformation principle. As a result they were forced to use essentially the
topological definition (\ref{a3.2}) of the curve, and the distance geometry
ceases to be a purely distance geometry, because it used non-metrical
concepts.

T-geometry uses the deformation principle essentially, and it does not need
the constraints of the type (\ref{a1.3}), (\ref{a1.4}). At first, in
T-geometry the deformation principle was applied to the $n$-dimensional
planes $(n=1,2,...)$ \cite{R2000}. In the Euclidean geometry the $n$%
-dimensional plane may be introduced as a linear envelope of $n$ linear
independent vectors. As far as the condition of the linear dependence in the
Euclidean geometry may be written in the form (\ref{a3.5}), the $n$%
-dimensional plane $\mathcal{T}_{\mathcal{P}^{n}}$, determined by $n+1$
points $\mathcal{P}^{n}=\left\{ P_{0},P_{1},...P_{n}\right\} $, may be
represented as the set of points 
\begin{equation}
\mathcal{T}_{\mathcal{P}^{n}}=\mathcal{T}\left( \mathcal{P}^{n}\right)
=\left\{ R|F_{n+1}\left( \mathcal{P}^{n},R\right) =0\right\} ,\qquad
F_{n}\left( \mathcal{P}^{n}\right) \neq 0  \label{a3.16}
\end{equation}
where $F_{n}$ is the Gram determinant, defined by the relation (\ref{a3.5}).
As far as the definition (\ref{a3.16}) contains only the world function, the
relation (\ref{a3.16}) may be expanded on the arbitrary $\sigma $-space $%
V=\left\{ \sigma ,\Omega \right\} $, and it is an application of the
deformation principle to the Euclidean planes. Further the action of the
deformation principle was extended to all geometrical objects \cite{R01},
and it became clear that the topological concept of the curve (\ref{a3.2})
is redundant at the geometry construction.

Application of the topological definition of the curve (\ref{a3.2}) leads
with necessity to the supposition on transitivity of the parallelism
property. In the nondegenerate T-geometry the property of parallelism is
intransitive, in general. Considering the Riemannian geometry as a special
case of the T-geometry, we may approach only local transitivity of the
parallelism. In other words, the property of the parallelism is transitive
only for vectors, having a common origin. The transitivity of the
parallelism cannot be achieved for the vectors at different points of the
space. To avoid the incompatibility between the demand of the transitivity,
which follows from the definition of the curve (\ref{a3.2}) and
impossibility of realizing this property for vectors at different points of
the space, the conventional Riemannian geometry is forced to reject the
absolute parallelism of vectors, remaining only parallelism at the same
point, where the transitivity of parallelism is possible.

\section{Concluding remarks}

Thus, the crisis in the science is a social-scientific phenomenon, where
both the social aspect and the scientific one are equally important. The
origin of the crisis is a mistake in the fundamental statements of a theory,
which could not be discovered for a long time. As a result the mistake
develops into a preconception. But the mistake as well as the preconception
is not yet a crisis. Mistakes are natural products of the human activity.
They are to be discovered and corrected. The crises appears, if the
scientific community does not want to accept the new version of the theory,
which appears after correction of the mistake (or preconception). Usually
the long existence of the erroneous theory is connected with the existence
of the Ptolemaic approach to the theory and with the Ptolemaic criterion of
this theory evaluation. The Newtonian approach does not yeasay the long
existence of an erroneous theory.

The described crisis in the geometry is not the first crisis and,
apparently, it is not the last one. In such a situation some recommendation,
how to overcome the crisis in the existing conception might be useful. My
advice for the scientific community is rather simple. If we have the crisis
phenomena in the fundamental theory, and there exist a new theory pretending
to correction of mistakes, the new theory should be evaluated basing on the
Newtonian criterion (but not on the Ptolemaic one). It should keep in mind
also the simple circumstance, that some corollaries of the new theory may
depreciate results, obtained on the basis of the existing erroneous theory.
(In the given case some papers dealing with foundation of the geometry on
the basis of the topological concepts may be depreciated). Authors of such
papers are not inclined to accept the new conception.

Practically I cannot give advices for researchers who will risk to find
future mistakes in the fundamental statements, because I do not know, how I
succeeded to find mistakes and to overcome the existing preconceptions. At
any rate, my ideas on, why I succeeded to find mistakes and to construct the
T-geometry, are rather cloudy and subjective. Nevertheless, I shall try to
formulate the principal circumstance, which allow one to construct the
T-geometry.

At first, I followed the Newtonian investigation strategy: "Hypotheses non
fingo" and evaluated my results according to the Newtonian criteria. Second,
I am physicist and I appear to be free of mathematical and geometrical
preconceptions, although, apparently, I have my own physical preconceptions.
Besides, I considered the geometry to be only a tool for investigation of
the space-time and of the physical phenomena in the space-time. The last
circumstance helped me to overcome the negative relation of professional
geometers to my papers on geometry. (It was more difficult to overcome the
negative relations of physicists to my physical papers). Finally, I tried to
consider objections of my opponents (which were inequitable in many cases)
as an display of the objective cognition process at the stage of crisis. It
helped me to perceive existence of objections of my opponents as an
objective necessity, generated by the crisis.

\[
\]


\begin{thebibliography}{99}
\bibitem{N2000} S.P. Novikov, Second half of XXth century and its result:
crisis of the physics-mathematics community. \textit{The Transactions of the
Institute for the History of the Natural Sciences and Technology}, 2002. (in
Russian).(English translation made by A.Sossinsky \textit{AMS Translations},
S.Novikov's Seminar, 2004.

\bibitem{R2004} Yu. A. Rylov, Coordinateless description and deformation
principle as a foundations of physical geometry. (Available at
http://arXiv.org/ abs/math.GM/ 0312160)

\bibitem{R03} Yu. A. Rylov, Model conception of quantum phenomena: logical
structure and investigation methods. (In preparation, available at
http://arXiv.org/abs /physics/0310050, v2).

\bibitem{R62} Yu.A. Rylov, On a possibility of the Riemannian space
description in terms of a finite interval. \textit{Izvestiya Vysshikh
Uchebnych Zavedenii, Matematika}. No.3(28), 131-142, (1962). (in Russian)

\bibitem{R64} Yu.A.Rylov, : Relative gravitational field and conservation
laws in general relativity. \textit{Ann. Phys.}(Leipzig) \textbf{12},
329-353, (1964).

\bibitem{S60} J. L.~Synge, \textit{Relativity: The General Theory},
North-Holland, Amsterdam, 1960.

\bibitem{R90} Yu. A.~Rylov, Extremal properties of Synge's world function
and discrete geometry. \textit{J. Math. Phys.} \textbf{31}, 2876-2890,
(1990).

\bibitem{R91} Yu. A.~Rylov, Non-Riemannian model of space-time responsible
for quantum effects. \textit{J. Math. Phys.} \textbf{32}, 2092-2098, (1991).

\bibitem{R92} Yu.A. Rylov, "Distorted Riemannian space and technique of
differential geometry". \textit{J.Math. Phys}. \textbf{33}, 4220-4224, (1992)

\bibitem{M28} K.~Menger, Untersuchen \"{u}ber allgemeine Metrik, \textit{%
Mathematische Annalen,} \textbf{100}, 75-113, (1928).

\bibitem{B53} L. M.~Blumenthal, \textit{Theory and Applications of Distance
Geometry}, Oxford, Clarendon Press, 1953.

\bibitem{R2000} Yu.A. Rylov Metric space: classification of finite subspaces
instead of constraints on metric. \textit{Proceedings on analysis and
geometry}, Novosibirsk, Publishing House of Mathematical institute, 2000.
pp. 481-504, (in Russian). English version: (Available at
http://arXiv.org/abs/math.MG/9905111).

\bibitem{R01} Yu.A. Rylov, Geometry without topology as a new conception of
geometry.\textit{\ Int. Jour. Mat. \& Mat. Sci.} \textbf{30}, iss. 12,
733-760, (2002), (available at http://arXiv.org/abs/math.MG/0103002).
\end{thebibliography}
\end{document}